# A Novel Approach for Flexible Body Dynamics Computation via Synthesizing Incremental Motions in Reconfigured Inertial Frames

Xiaobo Liu

Sandia National Laboratories

Albuquerque, NM 87123

e-mail: xialiu@sandia.gov

## Abstract

A novel approach is presented for computing flexible body dynamics based on conventional structural dynamics models. This approach innovatively captures the rigid body motion component embedded within a flexible body's movement, generates and synthesizes single-step responses in a sequence of reconfigured inertial frames that follow the rigid body motion. By doing so, it effectively bypasses the complexities associated with modeling flexibilities and formulating highly nonlinear coupled motion equations. In addition to improving predictive accuracy, this approach offers valuable insights into the interaction between rigid body motion and structural vibration. By bridging these two aspects, it advances the understanding of flexible body dynamics and delivers a precise, efficient simulation framework for a wide range of engineering applications.

## 1. Introduction

The significance of flexible body dynamics research spans multiple application domains, including aerospace engineering, automotive design, robotics, and biomechanics. In these applications, materials and structures undergo coupled rigid body motion and structural vibration, and they must be designed to endure dynamic loads under a variety of operating conditions. The understanding of flexible dynamics is essential for ensuring safety, performance and reliability. Unlike rigid body dynamics, where the shape and size of a body remain unchanged throughout its motion, flexible body dynamics involves an intricate interplay between rigid body motion and vibrational modes of the structure. Despite the advancements made in this area, researchers continue to face challenges in achieving numerical techniques that are both accurate and computationally efficient. Two prominent methodologies to incorporate flexibility into large overall rigid body motion are the Floating Frame of Reference Formulation (FFRF) [1-5] and the Absolute Nodal Coordinate Formulation (ANCF) [6-11]. Both approaches are grounded in Lagrangian dynamics and are extensions of the finite element method (FEM). They utilize many core techniques from FEM, including spatial discretization, interpolation through shape functions, formulation of energy terms, derivation of mass and stiffness matrices. The complexities and accuracies of these formulations are largely influenced by the selections of generalized coordinates, shape functions and discretization schemes. Consequently, the choice of method often hinges on the specific structure, application and context of analysis, necessitating trade-offs between accuracy, modeling complexity and computational efficiency. As the field continues to evolve, ongoing research aims to address the persistent challenges, refine the existing methods to broaden their applicability and capabilities.

Interestingly, no method has been proposed in literature to analyze flexible body dynamics directly using conventional structural dynamics models derived from finite element method (FEM) or other discretization techniques. Although conventional models are widely used in structural vibration analysis, being derived from the small-strain assumption limits their validity to configurations near the undeformed state. As a result, they are generally applicable only to small vibrational motions and are inadequate for predicting long-term motions involving large rigid body displacements. A conventional modal analysis does generate rigid body modes, but these rigid body modes are typically only utilized for assessing model integrity and are excluded from the computation of dynamic responses. The reason is that rigid body modes represent only the patterns of rigid body deformation, or more precisely, they serve as the base vectors of rigid body motion velocities. Currently there is no established procedure to obtain spatial rigid body translation and rotation from these modes. To overcome these limitations, in this paper we propose an innovative approach that can effectively compute flexible body dynamics within the existing structural dynamics framework, enabling a more in-depth understanding of the coupled motion. To address the limitation imposed by the small-strain assumption, we propose a technique that employs a series of reconfigured inertial frames, which is inspired by a simple technique used in rectilinear motion analysis of a mass particle as explained below.

Suppose a particle $m$ is subjected to a constant force $f$, and is moving on a straight line with initial displacement $x_0$ and initial velocity $v_0$. The acceleration in the inertial coordinate system $OX$ is $a = f/m$. Let the time step be $\Delta t$. The velocities and displacements of the particle at adjacent time points are related by the following equations

$$v_{n+1} = v_n + a\Delta t, \tag{1}$$

$$x_{n+1} = x_n + \Delta x_n, \tag{2}$$

in which

$$\Delta x_n = v_n \Delta t + \frac{1}{2} a (\Delta t)^2. \tag{3}$$

We can use an alternative technique to analyze the motion. At time $t_n$, the particle has reached a position in $OX$, we can define a new inertial frame $OX|_n$ with its origin at the particle's position at this time. In this shifted frame, the equation of motion $f = ma$ remains unchanged, the particle's velocity $v_n$ serves as the new initial condition to predict its next state at $t_{n+1}$. The velocity $v_{n+1}$ is computed using Eq. (1), and the new displacement $\Delta x_n$ is obtained from Eq. (3). This process can be iteratively applied across the full time span. The global motion is then reconstructed by chaining together the single-step predictions in the successively shifted inertial frames.

While an iterative procedure is unnecessary for the simple point mass in rectilinear motion, it has motivated the development of a promising strategy addressing the challenges and limitations associated with complex models. By identifying a suitable sequence of inertial frames, it becomes possible to perform single-step predictions within these frames using the same model without violating assumptions on model validity, instead of attempting an inaccurate long-term prediction in a fixed inertial frame. The procedure begins with a single-step prediction in the original inertial frame, followed by an innovative method to extract the rigid body motion component embedded within the flexible body's movement. From this extracted rigid body motion, we can calculate the incremental rotation from linear motion quantities using the formulations derived in [12]. The inertial frame is then reconfigured to follow the rigid body motion through incorporation of both rigid body displacement and reorientation. Since the equations of motion remain valid in the reconfigured frame, the same model can be used for the next single-step prediction with new initial conditions. This process is repeated throughout the simulation time span, and the global motion is reconstructed by synthesizing the individual single-step results. This approach extends the applicability of conventional structural dynamics models, it effectively circumvents the complexities in modeling flexible bodies and formulating highly nonlinear coupled motion equations. In addition to enhancing predictive accuracy, it provides deeper insight into the interaction between rigid body motion and structural vibration. By bridging these two aspects, the proposed approach advances the understanding of flexible body dynamics, offers a robust and efficient simulation framework for a wide range of engineering applications.

The balance of this article is organized as follows. In Section 2, we review the main results on computation of rotation from linear motion quantities [12]; In Section 3, we give detailed derivation of the new method for flexible body dynamic analysis; An example is provided in Section 4 to validate the proposed method. Finally, concluding remarks are given in Section 5.

## 2. Computation of rotation from linear motion quantities

In this section we summarize the main results developed in [12] on the computation of rigid body rotation from linear motion quantities.

A rigid body is shown in Fig. 1. Let $OX_1X_2X_3$ be the global inertial frame with base unit vectors $\mathbf{n}_1, \mathbf{n}_2, \mathbf{n}_3$, $O'X'_1X'_2X'_3$ be a body-fixed frame with base unit vectors $\mathbf{n}'_1, \mathbf{n}'_2, \mathbf{n}'_3$. Three non-colinear points $P, Q, R$ are chosen on the body with $R$ being the reference point. The relative position vectors are formed as $\mathbf{p} = \overrightarrow{RP}$ and $\mathbf{q} = \overrightarrow{RQ}$. The magnitudes of the vectors are $p = \|\mathbf{p}\|$, $q = \|\mathbf{q}\|$.

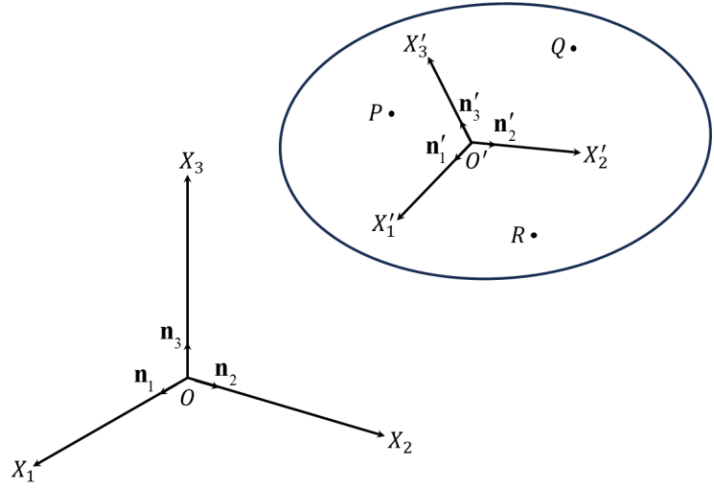

Fig. 1 A rigid body

Define unit vectors $\mathbf{e}_1 = \dfrac{\mathbf{p}}{p}$, $\mathbf{e}_2 = \dfrac{\mathbf{q}}{q}$, $\mathbf{e}_3 = \dfrac{\mathbf{e}_1 \times \mathbf{e}_2}{\|\mathbf{e}_1 \times \mathbf{e}_2\|}$. Let $\mathbf{B}$ be the transformation matrix such that

$$\begin{Bmatrix} \mathbf{e}_1 \\ \mathbf{e}_2 \\ \mathbf{e}_3 \end{Bmatrix} = \mathbf{B} \begin{Bmatrix} \mathbf{n}'_1 \\ \mathbf{n}'_2 \\ \mathbf{n}'_3 \end{Bmatrix}. \tag{4}$$

The rows of matrix $\mathbf{B}$ are components of the vectors $\mathbf{e}_1, \mathbf{e}_2, \mathbf{e}_3$ in the body-fixed frame. Three frames have been introduced: the global frame $OX_1 X_2 X_3$, the body-fixed frame $O'X'_1 X'_2 X'_3$, and the generally oblique frame spanned by $\mathbf{e}_1, \mathbf{e}_2, \mathbf{e}_3$. To be consistent in notations, for any vector $\mathbf{v}$, we use $(v^1, v^2, v^3)$ for its components in the global frame, $(v'^1, v'^2, v'^3)$ for its components in the body-fixed frame, and $(\hat{v}^1, \hat{v}^2, \hat{v}^3)$ for its components in the frame spanned by $\mathbf{e}_1, \mathbf{e}_2, \mathbf{e}_3$. Let $\mathbf{v}_{P/R}$ be the relative velocity of $P$ to $R$, $\mathbf{v}_{Q/R}$ be the relative velocity of $Q$ to $R$, $\mathbf{a}_{P/R}$ be the relative acceleration of $P$ to $R$, and $\mathbf{a}_{Q/R}$ be the relative acceleration of $Q$ to $R$. Then the angular velocity components can be calculated as

$$\hat{\omega}^3 = \dfrac{1}{p\|\mathbf{e}_1 \times \mathbf{e}_2\|} \mathbf{v}_{P/R} \cdot \left( \mathbf{e}_2 - \mathbf{e}_1 (\mathbf{e}_1 \cdot \mathbf{e}_2) \right), \tag{5}$$

$$\hat{\omega}^2 = \dfrac{-\mathbf{v}_{P/R} \cdot \mathbf{e}_3}{p\|\mathbf{e}_1 \times \mathbf{e}_2\|}, \tag{6}$$

$$\hat{\omega}^1 = \frac{\mathbf{v}_{Q/R} \cdot \mathbf{e}_3}{q \|\mathbf{e}_1 \times \mathbf{e}_2\|}, \tag{7}$$

$$\begin{Bmatrix} \omega'^1 \\ \omega'^2 \\ \omega'^3 \end{Bmatrix} = \mathbf{B}^T \begin{Bmatrix} \hat{\omega}^1 \\ \hat{\omega}^2 \\ \hat{\omega}^3 \end{Bmatrix}. \tag{8}$$

To simplify notation, we introduce a function that encapsulates Eqs. (5)-(8),

$$\boldsymbol{\omega} = g\left(\mathbf{v}_{P/R}, \mathbf{v}_{Q/R}\right). \tag{9}$$

The angular acceleration components are

$$\hat{\alpha}^3 = \frac{\left(\mathbf{e}_3 \times (\boldsymbol{\omega} \times \mathbf{v}_{P/R}) - \mathbf{e}_3 \times \mathbf{a}_{P/R}\right) \cdot \mathbf{e}_1}{p}, \tag{10}$$

$$\hat{\alpha}^2 = \frac{\left(\boldsymbol{\omega} \times \mathbf{v}_{P/R} - \mathbf{a}_{P/R}\right) \cdot \mathbf{e}_3}{p \|\mathbf{e}_1 \times \mathbf{e}_2\|}, \tag{11}$$

$$\hat{\alpha}^1 = \frac{\left(\mathbf{a}_{Q/R} - \boldsymbol{\omega} \times \mathbf{v}_{Q/R}\right) \cdot \mathbf{e}_3}{q \|\mathbf{e}_1 \times \mathbf{e}_2\|}, \tag{12}$$

$$\begin{Bmatrix} \alpha'^1 \\ \alpha'^2 \\ \alpha'^3 \end{Bmatrix} = \mathbf{B}^T \begin{Bmatrix} \hat{\alpha}^1 \\ \hat{\alpha}^2 \\ \hat{\alpha}^3 \end{Bmatrix}. \tag{13}$$

Similarly, we use a function to encapsulate Eqs. (10)-(13),

$$\boldsymbol{\alpha} = h\left(\mathbf{v}_{P/R}, \mathbf{v}_{Q/R}, \mathbf{a}_{P/R}, \mathbf{a}_{Q/R}, \boldsymbol{\omega}\right). \tag{14}$$

The evolutions of relative velocity components in the body-fixed frame are governed by the following differential equations,

$$\begin{Bmatrix} \dot{v}'^1_{P/R} \\ \dot{v}'^2_{P/R} \\ \dot{v}'^3_{P/R} \end{Bmatrix} = \begin{Bmatrix} a'^1_{P/R} \\ a'^2_{P/R} \\ a'^3_{P/R} \end{Bmatrix} - \begin{Bmatrix} \omega'^2 v'^3_{P/R} - \omega'^3 v'^2_{P/R} \\ \omega'^3 v'^1_{P/R} - \omega'^1 v'^3_{P/R} \\ \omega'^1 v'^2_{P/R} - \omega'^2 v'^1_{P/R} \end{Bmatrix}, \tag{15}$$

$$\begin{Bmatrix} \dot{v}'^1_{Q/R} \\ \dot{v}'^2_{Q/R} \\ \dot{v}'^3_{Q/R} \end{Bmatrix} = \begin{Bmatrix} a'^1_{Q/R} \\ a'^2_{Q/R} \\ a'^3_{Q/R} \end{Bmatrix} - \begin{Bmatrix} \omega'^2 v'^3_{Q/R} - \omega'^3 v'^2_{Q/R} \\ \omega'^3 v'^1_{Q/R} - \omega'^1 v'^3_{Q/R} \\ \omega'^1 v'^2_{Q/R} - \omega'^2 v'^1_{Q/R} \end{Bmatrix}. \tag{16}$$

The orientation of the rigid body is described by the coordinate transformation matrix of the body-fixed frame $O'X'_1X'_2X'_3$ with respect to the global frame $OX_1X_2X_3$. Any change in the configuration of a rigid body can be achieved by a rotation about an appropriate direction through an appropriate angle. Let $\boldsymbol{\lambda} = \begin{bmatrix} \lambda^1 & \lambda^2 & \lambda^3 \end{bmatrix}^T$ be the unit vector in the direction of rotation with components expressed in the global frame, $\theta$ be the rotation angle, and let

$$\hat{\mathbf{M}} = \begin{bmatrix} 0 & \lambda^3 & -\lambda^2 \\ -\lambda^3 & 0 & \lambda^1 \\ \lambda^2 & -\lambda^1 & 0 \end{bmatrix}, \tag{17}$$

then the coordinate transformation matrix is given by

$$\mathbf{A} = e^{\theta \hat{\mathbf{M}}}. \tag{18}$$

The matrix exponential can be calculated as

$$\begin{aligned} e^{\theta \hat{\mathbf{M}}} &= \mathbf{I} + \theta \hat{\mathbf{M}} + \frac{1}{2!}\theta^2 \hat{\mathbf{M}}^2 - \frac{1}{3!}\theta^3 \hat{\mathbf{M}} - \frac{1}{4!}\theta^4 \hat{\mathbf{M}}^2 + \frac{1}{5!}\theta^5 \hat{\mathbf{M}} + \frac{1}{6!}\theta^6 \hat{\mathbf{M}}^2 + \cdots \\ &= \mathbf{I} + \hat{\mathbf{M}}\left(\theta - \frac{1}{3!}\theta^3 + \frac{1}{5!}\theta^5 - \cdots\right) + \hat{\mathbf{M}}^2\left(\frac{1}{2!}\theta^2 - \frac{1}{4!}\theta^4 + \frac{1}{6!}\theta^6\right) \\ &= \mathbf{I} + \sin\theta \hat{\mathbf{M}} + (1-\cos\theta)\hat{\mathbf{M}}^2. \end{aligned} \tag{19}$$

Let $\mathbf{M} = \theta \hat{\mathbf{M}}$. The evolution of $\mathbf{M}$ is governed by the following equation

$$\dot{\mathbf{M}} = \left(f(\theta, \hat{\mathbf{M}})\right)^{-1} \boldsymbol{\Omega}, \tag{20}$$

where

$$\boldsymbol{\Omega} = \begin{bmatrix} 0 & \omega'^3 & -\omega'^2 \\ -\omega'^3 & 0 & \omega'^1 \\ \omega'^2 & -\omega'^1 & 0 \end{bmatrix}. \tag{21}$$

and

$$\left(f\left(\theta,\hat{\mathbf{M}}\right)\right)^{-1}=\begin{cases}\mathbf{I}, & \text{when } \theta=0; \\ \mathbf{I}-(\theta/2)ad_{\hat{\mathbf{M}}}-\left(\dfrac{\theta\cot(\theta/2)}{2}-1\right)\left(ad_{\hat{\mathbf{M}}}\right)^{2}, & \text{when } \theta\neq 0.\end{cases} \qquad (22)$$

In Eq. (22), the adjoint operator is defined as $ad_{\hat{\mathbf{M}}}=-\left[\hat{\mathbf{M}},\bullet\right]$, in which the Lie bracket is given by $\left[\hat{\mathbf{M}},\hat{\mathbf{N}}\right]=\hat{\mathbf{N}}\cdot\hat{\mathbf{M}}-\hat{\mathbf{M}}\cdot\hat{\mathbf{N}}$.

## 3. Flexible body dynamics

Supposed we have a discretized model of flexible body. Let $O\underline{X}_1\underline{X}_2\underline{X}_3$ be the global analysis inertial frame in which the motion quantities are eventually sought after, and $OX_1X_2X_3$ be the inertial frame in which the flexible boy is discretized. $OX_1X_2X_3$ is stationary in $O\underline{X}_1\underline{X}_2\underline{X}_3$. A particle on the flexible body is uniquely identified by its fixed material coordinates in $OX_1X_2X_3$, i.e., coordinates of the particle when the body is undeformed. We also introduce a dummy rigid body consisting of the massless material points in undeformed configuration. Let a frame $O'X_1'X_2'X_3'$ be fixed on the dummy rigid body. At the starting point, $O'X_1'X_2'X_3'$ coincides with $OX_1X_2X_3$. Pick a reference point $R$ and points $P$, $Q$ on the body for angular velocity and angular acceleration computations.

A small motion of the flexible body in $OX_1X_2X_3$ is governed by a system of second order differential equations,

$$\mathbf{M}_{\mathbf{x}}\ddot{\mathbf{x}}+\mathbf{C}_{\mathbf{x}}\dot{\mathbf{x}}+\mathbf{K}_{\mathbf{x}}\mathbf{x}=\mathbf{f}_{\mathbf{x}}, \qquad (23)$$

where $\mathbf{x}$ is the physical displacement vector with respect to the undeformed state of the body, $\mathbf{M}_{\mathbf{x}}$ is the mass matrix, $\mathbf{K}_{\mathbf{x}}$ is the stiffness matrix, , $\mathbf{f}_{\mathbf{x}}$ is the force vector. $\mathbf{C}_{\mathbf{x}}$ is the damping matrix. Typically, in structural dynamics applications, structural damping or proportional damping is assumed, or modal damping is assigned to the modes of interest. Put together the mode shapes of interest with the first few columns being the rigid body modes determined by constraints and boundary conditions,

$$\mathbf{\Psi}=\begin{bmatrix}\mathbf{\psi}_1 & \mathbf{\psi}_2 & \cdots & \mathbf{\psi}_n\end{bmatrix}. \qquad (24)$$

Let $\Psi_R$ be the sub-matrix consisting of rigid body modes. The modal mass, modal stiffness and modal damping matrices are computed as

$$\mathbf{M}_q = \Psi^T \mathbf{M}_x \Psi, \tag{25}$$

$$\mathbf{K}_q = \Psi^T \mathbf{K}_x \Psi, \tag{26}$$

$$\mathbf{C}_q = \Psi^T \mathbf{C}_x \Psi. \tag{27}$$

Introduce vector $\mathbf{q}$ of modal coordinates and assume modal expansion

$$\mathbf{x} = \Psi \mathbf{q}, \tag{28}$$

$$\dot{\mathbf{x}} = \Psi \dot{\mathbf{q}}. \tag{29}$$

Conversely, the modal coordinates and modal velocities can be computed as

$$\mathbf{q} = \mathbf{M}_q^{-1} \Psi^T \mathbf{M}_x \mathbf{x}, \tag{30}$$

$$\dot{\mathbf{q}} = \mathbf{M}_q^{-1} \Psi^T \mathbf{M}_x \dot{\mathbf{x}}. \tag{31}$$

Eqs. (30) and (31) will hereafter be referred to as the 'inverse modal expansion'. Eq. (23) can be condensed into the modal space,

$$\mathbf{M}_q \ddot{\mathbf{q}} + \mathbf{C}_q \dot{\mathbf{q}} + \mathbf{K}_q \mathbf{q} = \mathbf{f}_q, \tag{32}$$

in which $\mathbf{f}_q = \Psi^T \mathbf{f}_x$. Define modal state vector

$$\boldsymbol{\eta} = \begin{Bmatrix} \mathbf{q} \\ \dot{\mathbf{q}} \end{Bmatrix}, \tag{33}$$

and let

$$\mathbf{H} = \begin{bmatrix} \mathbf{0} & \mathbf{I} \\ -\mathbf{M}_q^{-1}\mathbf{K}_q & -\mathbf{M}_q^{-1}\mathbf{C}_q \end{bmatrix}, \tag{34}$$

then Eqs. (32) can be converted into state-space form

$$\dot{\boldsymbol{\eta}} = \mathbf{H}\boldsymbol{\eta} + \begin{Bmatrix} \mathbf{0} \\ \mathbf{M}_q^{-1}\mathbf{f}_q \end{Bmatrix}. \tag{35}$$

We have obtained two evolution equations: Eq. (20) and Eq. (35). Similar to the shifted inertial frame technique introduced earlier for a single mass particle, we can use Eq. (20) to compute inertial frame evolution resulted from the embedded rigid body motion component of flexible body motion, and use Eq. (35) to obtain single-step predications in the reconfigured inertial frames throughout the entire time span, and finally synthesize the single-step predications to obtain the overall motion. This process involves decomposition of flexible body motion to obtain the embedded rigid body motion, and utilization of the formulations summarized in Section 2 to obtain rotation. The implementation steps are explained below.

Let the computation starts from $t_0$ with time step $\Delta t$. Suppose we are given conditions $\mathbf{x}_0$ and $\dot{\mathbf{x}}_0$, where the initial displacement $\mathbf{x}_0$ specifies the initial flexible deformation of the body. Calculate $\boldsymbol{\eta}_0$ using inverse modal expansion Eqs. (30) and (31), and $\dot{\boldsymbol{\eta}}_0$ from Eq. (35). Extract from $\dot{\boldsymbol{\eta}}_0$ the modal velocities and modal accelerations associated with rigid body modes, and use modal expansion to obtain the rigid body motion velocities and accelerations. Calculate relative velocities and relative accelerations, then use Eqs. (9) and (14) to obtain the initial angular velocity and angular acceleration of the dummy rigid body.

- Calculate one step change $\Delta\boldsymbol{\eta}_0$ in modal state vector from Eq. (35) using the fourth-Order Runge-Kutta method (RK4). Use modal expansion to obtain the change in physical state $\begin{Bmatrix} \Delta\mathbf{x}_0 \\ \Delta\dot{\mathbf{x}}_0 \end{Bmatrix}$. A subset of $\Delta\boldsymbol{\eta}_0$ is $\Delta\mathbf{q}_{R,0}$ containing the change of modal coordinates associated with the embedded rigid body motion. Calculate the embedded rigid body displacement $\Delta\mathbf{x}_{R,0}$ from $\Delta\mathbf{q}_{R,0}$ using modal expansion. $\Delta\mathbf{x}_{R,0}$ is also the displacement of the dummy rigid body,

- Apply the RK4 procedure to Eq. (20) to obtain $\mathbf{M}_1$ at $t_1$, and use Eq. (19) to obtain the orientation of the dummy rigid body at $t_1$. The displaced dummy rigid body with new orientation represents the undeformed state of flexible body at $t_1$. Let the reconfigured inertial frame fixed on the dummy body be denoted by $O'X'_1X'_2X'_3|_{t_1}$. Calculate the coordinate transformation matrix $\mathbf{R}_{1/0}$ of $O'X'_1X'_2X'_3|_{t_1}$ with respect to $O'X'_1X'_2X'_3$.

- Define a new inertial frame $OX_1X_2X_3|_{t_1}$ coinciding with $O'X'_1X'_2X'_3|_{t_1}$ at point $t_1$. In this new inertial frame, the governing equations of motion stay the same, the state vector is updated to $\begin{Bmatrix} \mathbf{x}_1 \\ \dot{\mathbf{x}}_1 \end{Bmatrix} = T_{1/0}\left(\begin{Bmatrix} \mathbf{x}_0 + \Delta\mathbf{x}_0 - \Delta\mathbf{x}_{R,0} \\ \dot{\mathbf{x}}_0 + \Delta\dot{\mathbf{x}}_0 \end{Bmatrix}\right)$, in which $T_{1/0}$ represents the coordinate transformation operation from $t_0$ to $t_1$. Depending on how the degrees of freedom (DOFs) are ordered in the discretized model, the transformation operation may involve reordering of the DOFs and multiplications by the coordinate transformation matrix $\mathbf{R}_{1/0}$. In the updated state vector, rigid body displacements have been removed because they are already accounted for in the reconfigured new inertial frame. Coordinate transformations has also been performed to bring the new state vector to the current inertial frame $OX_1X_2X_3|_{t_1}$.

- Use the updated state vector as new initial condition and follow the same steps, calculate the state at the next step $t_2 = t_1 + \Delta t$.

- Iterate to the end point.

The detailed RK4 procedure from $t_n$ to $t_{n+1}$ is given below. In the following context, $T_{m/n}$ represents the coordinate transformation from time $t_n$ to time $t_m$.

- At step $n$, suppose we have obtained the modal state vector $\mathbf{\eta}_n$, the orientation $\mathbf{M}_n$ of frame $O'X'_1X'_2X'_3|_n$. We can derived from $\mathbf{\eta}_n$ the relative velocities $(\mathbf{v}_{P/R,n}, \mathbf{v}_{Q/R,n})$ and relative accelerations $(\mathbf{a}_{P/R,n}, \mathbf{a}_{Q/R,n})$ associated with the embedded rigid body motion, with components expressed in $O'X'_1X'_2X'_3|_n$. The angular velocity and angular acceleration can be obtained using the functions as defined in Eqs. (9) and (14),

$$\mathbf{\omega}_n = g\left(\mathbf{v}_{P/R,n}, \mathbf{v}_{Q/R,n}\right), \tag{36}$$

$$\mathbf{\alpha}_n = h\left(\mathbf{v}_{P/R,n}, \mathbf{v}_{Q/R,n}, \mathbf{a}_{P/R,n}, \mathbf{a}_{Q/R,n}, \mathbf{\omega}_n\right). \tag{37}$$

- Step $k_1$,

$$\left.\begin{array}{c} \boldsymbol{\omega}_n, \mathbf{M}_n, \mathbf{v}_{P/R,n}, \mathbf{v}_{Q/R,n}, \mathbf{a}_{P/R,n}, \mathbf{a}_{Q/R,n} \\ \dot{\boldsymbol{\eta}}_n = \mathbf{H}\boldsymbol{\eta}_n + \begin{Bmatrix} \mathbf{0} \\ \mathbf{M}_\mathbf{q}^{-1}\mathbf{f}_{\mathbf{q},n} \end{Bmatrix} \end{array}\right\} \Rightarrow \begin{cases} \Delta_{k_1}\boldsymbol{\eta} = \dot{\boldsymbol{\eta}}_n \Delta t \text{ and subset } \Delta_{k_1}\mathbf{q}_R \\[1em] \Delta_{k_1}\mathbf{v}_{P/R} = \Delta t \left( \begin{Bmatrix} a'^{1}_{P/R,n} \\ a'^{2}_{P/R,n} \\ a'^{3}_{P/R,n} \end{Bmatrix} - \begin{Bmatrix} \omega'^{2}_n v'^{3}_{P/R,n} - \omega'^{3}_n v'^{2}_{P/R,n} \\ \omega'^{3}_n v'^{1}_{P/R,n} - \omega'^{1}_n v'^{3}_{P/R,n} \\ \omega'^{1}_n v'^{2}_{P/R,n} - \omega'^{2}_n v'^{1}_{P/R,n} \end{Bmatrix} \right) \\[2em] \Delta_{k_1}\mathbf{v}_{Q/R} = \Delta t \left( \begin{Bmatrix} a'^{1}_{Q/R,n} \\ a'^{2}_{Q/R,n} \\ a'^{3}_{Q/R,n} \end{Bmatrix} - \begin{Bmatrix} \omega'^{2}_n v'^{3}_{Q/R,n} - \omega'^{3}_n v'^{2}_{Q/R,n} \\ \omega'^{3}_n v'^{1}_{Q/R,n} - \omega'^{1}_n v'^{3}_{Q/R,n} \\ \omega'^{1}_n v'^{2}_{Q/R,n} - \omega'^{2}_n v'^{1}_{Q/R,n} \end{Bmatrix} \right) \\[2em] \Delta_{k_1}\mathbf{M} = \Delta t \left( \left( f\left(\theta_n, \hat{\mathbf{M}}_n\right) \right)^{-1} \boldsymbol{\Omega}_n \right) \end{cases}$$

- Step $k_2$,

$$\mathbf{v}_{P/R,k_2} = \mathbf{v}_{P/R,n} + \frac{1}{2}\Delta_{k_1}\mathbf{v}_{P/R}, \quad \mathbf{v}_{Q/R,k_2} = \mathbf{v}_{Q/R,n} + \frac{1}{2}\Delta_{k_1}\mathbf{v}_{Q/R}$$

$$\boldsymbol{\omega}_{k_2} = g\left(\mathbf{v}_{P/R,k_2}, \mathbf{v}_{Q/R,k_2}\right)$$

$$\mathbf{M}_{k_2} = \mathbf{M}_n + \frac{1}{2}\Delta_{k_1}\mathbf{M} \rightarrow \mathbf{R}_{k_2}$$

$$\left\{\begin{matrix}\mathbf{x}_{k_2}\\ \dot{\mathbf{x}}_{k_2}\end{matrix}\right\} = T_{k_2/n}\left(\left\{\begin{matrix}\mathbf{x}_n\\ \dot{\mathbf{x}}_n\end{matrix}\right\} + \begin{bmatrix}\boldsymbol{\Psi} & 0\\ 0 & \boldsymbol{\Psi}\end{bmatrix}\cdot\frac{1}{2}\Delta_{k_1}\boldsymbol{\eta} - \left\{\begin{matrix}\boldsymbol{\Psi}_R\cdot\frac{1}{2}\Delta_{k_1}\mathbf{q}_R\\ 0\end{matrix}\right\}\right)$$

Obtain $\boldsymbol{\eta}_{k_2}$ from $\left\{\begin{matrix}\mathbf{x}_{k_2}\\ \dot{\mathbf{x}}_{k_2}\end{matrix}\right\}$ using inverse modal expansion, calculate $\dot{\boldsymbol{\eta}}_{k_2} = \mathbf{H}\boldsymbol{\eta}_{k_2} + \left\{\begin{matrix}0\\ \mathbf{M}_\mathbf{q}^{-1}\left(\dfrac{\mathbf{f}_{\mathbf{q},n} + \mathbf{f}_{\mathbf{q},n+1}}{2}\right)\end{matrix}\right\}$

Use $\dot{\boldsymbol{\eta}}_{k_2}$ and modal expansion to obtain accelerations and relative accelerations $\left(\mathbf{a}_{P/R,k_2}, \mathbf{a}_{Q/R,k_2}\right)$

$$\Rightarrow \begin{cases} \Delta_{k_2}\boldsymbol{\eta} = \Delta t\dot{\boldsymbol{\eta}}_{k_2} \text{ and subset } \Delta_{k_2}\mathbf{q}_R \\ \Delta_{k_2}\mathbf{v}_{P/R} = \Delta t\left(\left\{\begin{matrix}a'^1_{P/R,k_2}\\ a'^2_{P/R,k_2}\\ a'^3_{P/R,k_2}\end{matrix}\right\} - \left\{\begin{matrix}\omega'^2_{k_2}v'^3_{P/R,k_2} - \omega'^3_{k_2}v'^2_{P/R,k_2}\\ \omega'^3_{k_2}v'^1_{P/R,k_2} - \omega'^1_{k_2}v'^3_{P/R,k_2}\\ \omega'^1_{k_2}v'^2_{P/R,k_2} - \omega'^2_{k_2}v'^1_{P/R,k_2}\end{matrix}\right\}\right) \\ \Delta_{k_2}\mathbf{v}_{Q/R} = \Delta t\left(\left\{\begin{matrix}a'^1_{Q/R,k_2}\\ a'^2_{Q/R,k_2}\\ a'^3_{Q/R,k_2}\end{matrix}\right\} - \left\{\begin{matrix}\omega'^2_{k_2}v'^3_{Q/R,k_2} - \omega'^3_{k_2}v'^2_{Q/R,k_2}\\ \omega'^3_{k_2}v'^1_{Q/R,k_2} - \omega'^1_{k_2}v'^3_{Q/R,k_2}\\ \omega'^1_{k_2}v'^2_{Q/R,k_2} - \omega'^2_{k_2}v'^1_{Q/R,k_2}\end{matrix}\right\}\right) \\ \Delta_{k_2}\mathbf{M} = \Delta t\left(\left(f\left(\theta_{k_2},\hat{\mathbf{M}}_{k_2}\right)\right)^{-1}\boldsymbol{\Omega}_{k_2}\right) \end{cases}$$

- Step $k_3$,

$$\left.\begin{aligned}
&\mathbf{v}_{P/R,k_3} = \mathbf{v}_{P/R,n} + \frac{1}{2}\Delta_{k_2}\mathbf{v}_{P/R}, \quad \mathbf{v}_{Q/R,k_3} = \mathbf{v}_{Q/R,n} + \frac{1}{2}\Delta_{k_2}\mathbf{v}_{Q/R} \\
&\boldsymbol{\omega}_{k_3} = g\left(\mathbf{v}_{P/R,k_3}, \mathbf{v}_{Q/R,k_3}\right) \\
&\mathbf{M}_{k_3} = \mathbf{M}_n + \frac{1}{2}\Delta_{k_2}\mathbf{M} \to \mathbf{R}_{k_3} \\
&\left\{\begin{matrix}\mathbf{x}_{k_3}\\ \dot{\mathbf{x}}_{k_3}\end{matrix}\right\} = T_{k_3/n}\left(\left\{\begin{matrix}\mathbf{x}_n\\ \dot{\mathbf{x}}_n\end{matrix}\right\} + \begin{bmatrix}\boldsymbol{\Psi} & 0 \\ 0 & \boldsymbol{\Psi}\end{bmatrix}\cdot\frac{1}{2}\Delta_{k_2}\boldsymbol{\eta} - \left\{\begin{matrix}\boldsymbol{\Psi}_R \cdot \frac{1}{2}\Delta_{k_2}\mathbf{q}_R \\ 0\end{matrix}\right\}\right) \\
&\text{Obtain } \boldsymbol{\eta}_{k_3} \text{ from } \left\{\begin{matrix}\mathbf{x}_{k_3}\\ \dot{\mathbf{x}}_{k_3}\end{matrix}\right\} \text{ using inverse modal expansion, calculate } \dot{\boldsymbol{\eta}}_{k_3} = \mathbf{H}\boldsymbol{\eta}_{k_3} + \left\{\begin{matrix}0 \\ \mathbf{M}_\mathbf{q}^{-1}\left(\dfrac{\mathbf{f}_{\mathbf{q},n}+\mathbf{f}_{\mathbf{q},n+1}}{2}\right)\end{matrix}\right\} \\
&\text{Use } \dot{\boldsymbol{\eta}}_{k_3} \text{ and modal expansion to obtain accelerations and relative accelerations } \mathbf{a}_{P/R,k_2}, \mathbf{a}_{Q/R,k_2}
\end{aligned}\right\}$$

$$\Rightarrow \begin{cases}
\Delta_{k_3}\boldsymbol{\eta} = \Delta t \dot{\boldsymbol{\eta}}_{k_3} \text{ and subset } \Delta_{k_3}\mathbf{q}_R \\
\Delta_{k_3}\mathbf{v}_{P/R} = \Delta t\left(\left\{\begin{matrix}a'^1_{P/R,k_3}\\ a'^2_{P/R,k_3}\\ a'^3_{P/R,k_3}\end{matrix}\right\} - \left\{\begin{matrix}\omega'^2_{k_3}v'^3_{P/R,k_3} - \omega'^3_{k_3}v'^2_{P/R,k_3}\\ \omega'^3_{k_2}v'^1_{P/R,k_3} - \omega'^1_{k_2}v'^3_{P/R,k_3}\\ \omega'^1_{k_2}v'^2_{P/R,k_3} - \omega'^2_{k_2}v'^1_{P/R,k_3}\end{matrix}\right\}\right) \\
\Delta_{k_3}\mathbf{v}_{Q/R} = \Delta t\left(\left\{\begin{matrix}a'^1_{Q/R,k_3}\\ a'^2_{Q/R,k_3}\\ a'^3_{Q/R,k_3}\end{matrix}\right\} - \left\{\begin{matrix}\omega'^2_{k_3}v'^3_{Q/R,k_3} - \omega'^3_{k_3}v'^2_{Q/R,k_3}\\ \omega'^3_{k_3}v'^1_{Q/R,k_3} - \omega'^1_{k_3}v'^3_{Q/R,k_3}\\ \omega'^1_{k_3}v'^2_{Q/R,k_3} - \omega'^2_{k_3}v'^1_{Q/R,k_3}\end{matrix}\right\}\right) \\
\Delta_{k_3}\mathbf{M} = \Delta t\left(\left(f\left(\boldsymbol{\theta}_{k_3}, \hat{\mathbf{M}}_{k_3}\right)\right)^{-1}\boldsymbol{\Omega}_{k_3}\right)
\end{cases}$$

- Step $k_4$,

$$\left.\begin{aligned}
&\mathbf{v}_{P/R,k_4} = \mathbf{v}_{n,P/R} + \Delta_{k_3}\mathbf{v}_{P/R}, \quad \mathbf{v}_{Q/R,k_4} = \mathbf{v}_{n,Q/R} + \Delta_{k_3}\mathbf{v}_{Q/R} \\
&\boldsymbol{\omega}_{k_4} = g\left(\mathbf{v}_{P/R,k_4}, \mathbf{v}_{Q/R,k_4}\right) \\
&\mathbf{M}_{k_4} = \mathbf{M}_n + \Delta_{k_3}\mathbf{M} \to \mathbf{R}_{k_4} \\
&\begin{Bmatrix}\mathbf{x}_{k_4}\\\dot{\mathbf{x}}_{k_4}\end{Bmatrix} = T_{k_4/n}\left(\begin{Bmatrix}\mathbf{x}_n\\\dot{\mathbf{x}}_n\end{Bmatrix} + \begin{bmatrix}\boldsymbol{\Psi} & \mathbf{0}\\\mathbf{0} & \boldsymbol{\Psi}\end{bmatrix}\cdot\Delta_{k_3}\boldsymbol{\eta} - \begin{Bmatrix}\boldsymbol{\Psi}_R\cdot\Delta_{k_3}\mathbf{q}_R\\\mathbf{0}\end{Bmatrix}\right) \\
&\text{Obtain } \boldsymbol{\eta}_{k_4} \text{ from } \begin{Bmatrix}\mathbf{x}_{k_4}\\\dot{\mathbf{x}}_{k_4}\end{Bmatrix} \text{ using inverse modal expansion, calculate } \dot{\boldsymbol{\eta}}_{k_4} = \mathbf{H}\boldsymbol{\eta}_{k_4} + \begin{Bmatrix}\mathbf{0}\\\mathbf{M}_\mathbf{q}^{-1}\mathbf{f}_{\mathbf{q},n+1}\end{Bmatrix}
\end{aligned}\right\}$$

$$\Rightarrow \begin{cases}\Delta_{k_4}\boldsymbol{\eta} = \Delta t\,\dot{\boldsymbol{\eta}}_{k_4}\\ \Delta_{k_4}\mathbf{M} = \Delta t\left(\left(f\left(\theta_{k_4},\hat{\mathbf{M}}_{k_4}\right)\right)^{-1}\boldsymbol{\Omega}_{k_4}\right)\end{cases}$$

Finally, orientation of the dummy body is updated to

$$\mathbf{M}_{n+1} = \mathbf{M}_n + \frac{1}{6}\left(\Delta_{k_1}\mathbf{M} + 2\Delta_{k_2}\mathbf{M} + 2\Delta_{k_3}\mathbf{M} + \Delta_{k_4}\mathbf{M}\right), \tag{38}$$

Then change in modal state vector is

$$\Delta_n\boldsymbol{\eta} = \frac{1}{6}\left(\Delta_{k_1}\boldsymbol{\eta} + 2\Delta_{k_2}\boldsymbol{\eta} + 2\Delta_{k_3}\boldsymbol{\eta} + \Delta_{k_4}\boldsymbol{\eta}\right). \tag{39}$$

From $\Delta_n\boldsymbol{\eta}$ we can obtain the subset $\Delta_n\mathbf{q}_R$ corresponding to the change in modal coordinates associated with the embedded rigid body motion. The physical state vector at $t_{n+1}$ is

$$\begin{Bmatrix}\mathbf{x}_{n+1}\\\dot{\mathbf{x}}_{n+1}\end{Bmatrix} = T_{n+1/n}\left(\begin{Bmatrix}\mathbf{x}_n\\\dot{\mathbf{x}}_n\end{Bmatrix} + \begin{bmatrix}\boldsymbol{\Psi} & \mathbf{0}\\\mathbf{0} & \boldsymbol{\Psi}\end{bmatrix}\cdot\Delta_n\boldsymbol{\eta} - \begin{Bmatrix}\boldsymbol{\Psi}_R\cdot\Delta_n\mathbf{q}_R\\\mathbf{0}\end{Bmatrix}\right). \tag{40}$$

Using inverse modal expansion we can obtain $\boldsymbol{\eta}_{n+1}$ from $\begin{Bmatrix}\mathbf{x}_{n+1}\\\dot{\mathbf{x}}_{n+1}\end{Bmatrix}$.

## 4. Numerical example

In this section, we present a numerical example to demonstrate the effectiveness of the proposed method. All modeling, programming, and numerical computations are performed in MATLAB, which imposes inherent limitations on matrix size. To accommodate these constraints, the simulation model is purposely kept at a level of complexity that is manageable within MATLAB, particularly in terms of nodal degrees of freedom. Nevertheless, the model remains sufficiently representative for validation and verification purposes. The simulation results demonstrate that the method accurately captures coupled rigid body motion and structural vibration.

A rectangular plate has dimensions $L_1 = 1000\ mm$, $L_2 = 50\ mm$, $L_3 = 10\ mm$. The material properties are: Young's modulus $E = 70000\ MPa$, Poisson's ratio $\upsilon = 0.33$, density $\rho = 2.7 \times 10^{-9}\ Mg/mm^3$. A concentrated mass of $0.001\ Mg$ is attached to each of the corner points: $P_1(-500, -25, 0)$, $P_2(500, -25, 0)$, $P_3(500, 25, 0)$, $P_4(-500, 25, 0)$, $P_5(-500, -25, 10)$, $P_6(500, -25, 10)$, $P_7(500, 25, 10)$, $P_8(-500, 25, 10)$. As shown in Fig. 2, the plate is meshed with tetra10 solid elements using MATLAB's PDE toolbox. The total mass and principal moments of inertia calculated from the finite element model are: $m = 0.0093\ Mg$, $I_{11} = 5.4839\ Mg \cdot mm^2$, $I_{22} = 2112.7024\ Mg \cdot mm^2$, $I_{33} = 2117.7732\ Mg \cdot mm^2$. The original mesh frame is chosen as the global analysis frame.

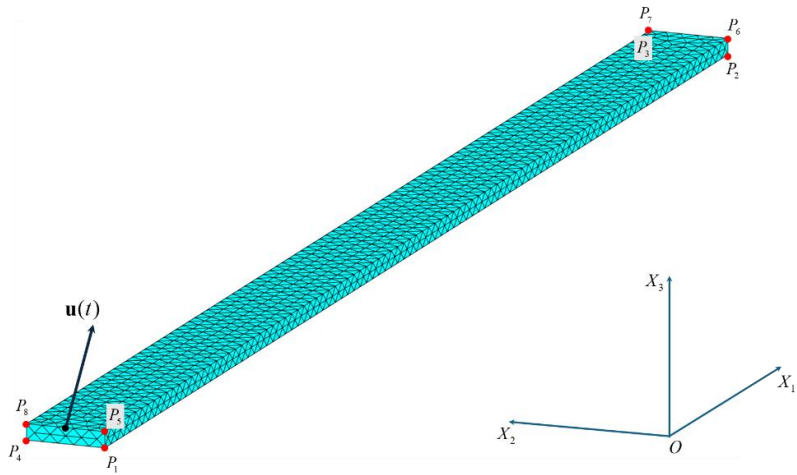

Fig. 2　FEA model of a rectangular flexible plate

The first six modes of the flexible plate are rigid body modes. The frequencies of the first ten flexible modes are given in Table 1, and the mode shapes are shown in Fig. 3. These flexible modes are used in computations together with the rigid body modes. We assume 5% modal damping for all the modes.

Table 1  Frequencies of the first ten flexible modes

| Mode # | Frequency (Hz) | Mode # | Frequency (Hz) |
|---|---|---|---|
| 7 | 24.5282 | 12 | 359.8918 |
| 8 | 85.8433 | 13 | 387.63 |
| 9 | 93.2656 | 14 | 534.3499 |
| 10 | 117.3671 | 15 | 544.946 |
| 11 | 206.1601 | 16 | 579.7381 |

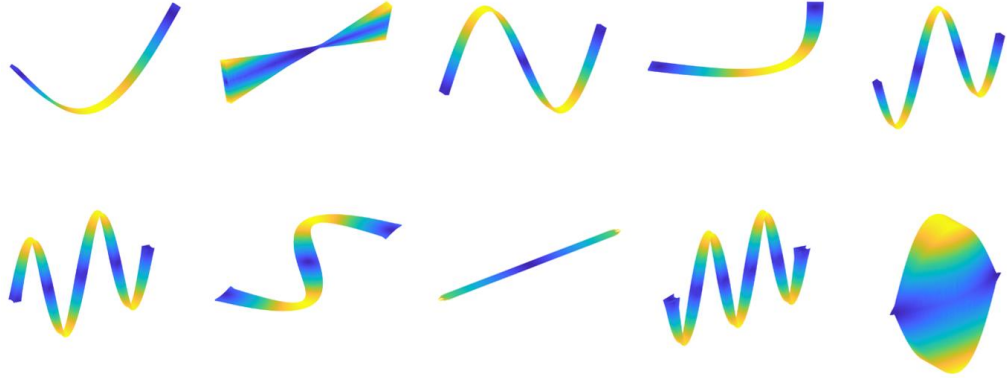

Fig. 3  Mode shapes of the first ten flexible modes

A $50\ N$ force is applied at the midpoint between $P_5$ and $P_8$ along direction $[1\ \ 0\ \ 1]'$ in the global analysis frame. Time history of the equivalent force components $u_1(t), u_3(t)$ is shown in Fig. 4. The dynamics of the flexible plate is computed using the method proposed in this paper with a time step $\Delta t = 10^{-4}$ second. $P_1$ is chosen as the reference point together with $P_2$ and $P_3$ for angular velocity and angular acceleration calculations. We capture the angular velocity and angular acceleration of the flexible plate as well as responses at selected points. For comparison, we also create a rigid plate model using the calculated mass and moment of inertia values. The same load is applied on the rigid plate. The motion of the rigid plate center of mass is governed by Newton's second law, the rotation of the rigid plate is governed by the angular momentum theorem. The component equations given by the angular moment theorem in body-fixed frame are the Euler equations. Numerical results are summarized in the remaining of this section.

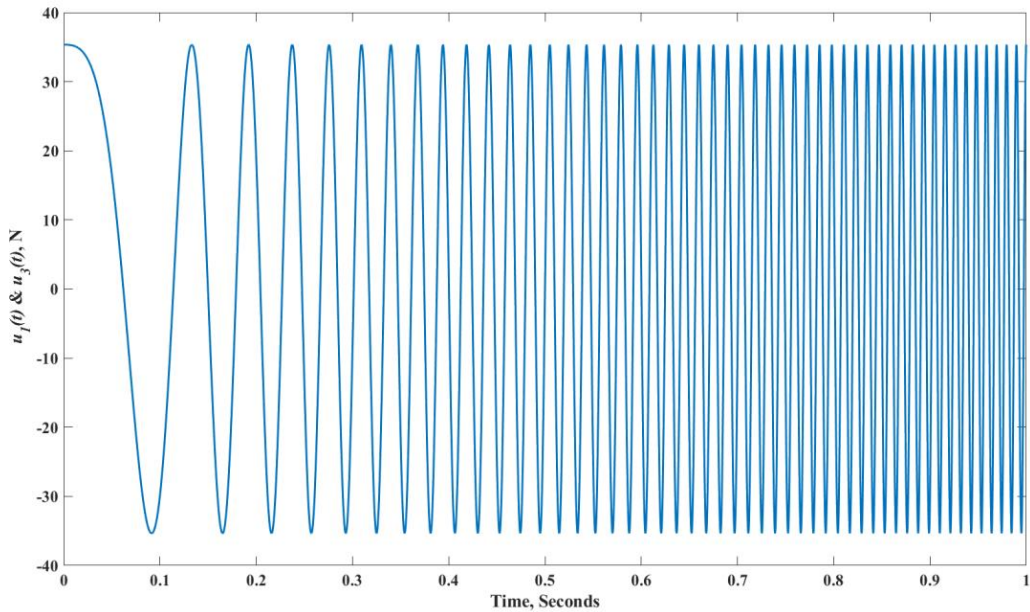

Fig. 4  Time history of applied force

For the rigid plate, the motion is translation in the $X_3$ direction with the center of mass, in conjunction with pure rotation about the body-fixed $X_2'$ axis. The embedded rigid body motion (RBM) component of the flexible plate resembles the rigid plate motion. As shown in Figs. 5 and 6, the angular acceleration and angular velocity components of the flexible plate in the body-fixed $X_1'$ and $X_3'$ directions are negligible. The displacements of the selected points $P_1$ through $P_8$ in the global $X_2$ direction are also trivial as shown in Fig. 7. As shown in Figs. 8 and 9, the angular acceleration and angular velocity components of the flexible plate in body-fixed $X_2'$ direction match precisely with those of the rigid plate.

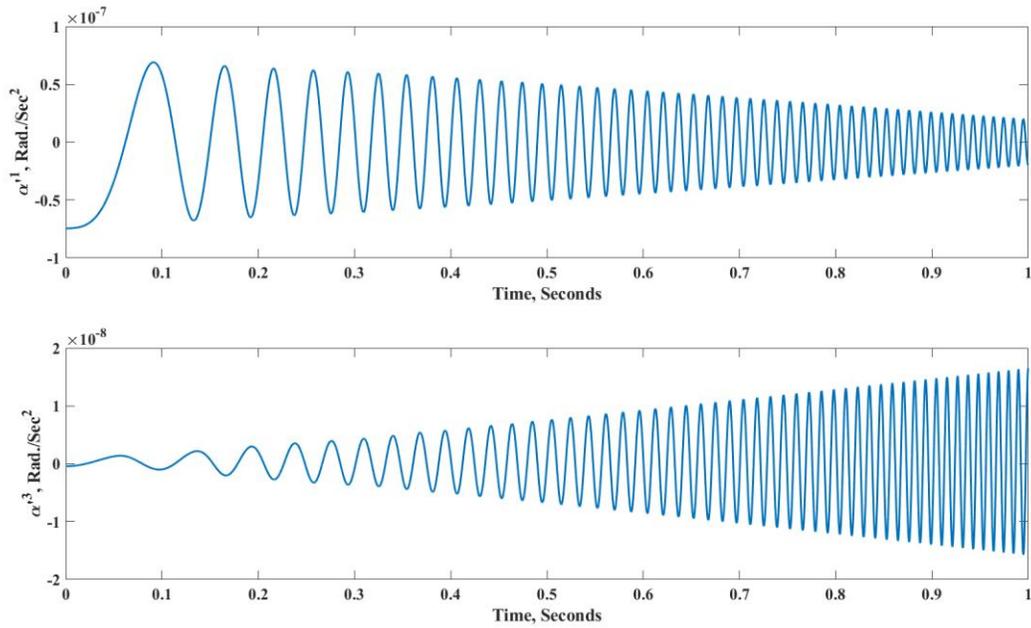

Fig. 5 Angular acceleration components of flexible plate in body-fixed $X'_1$ and $X'_3$ directions

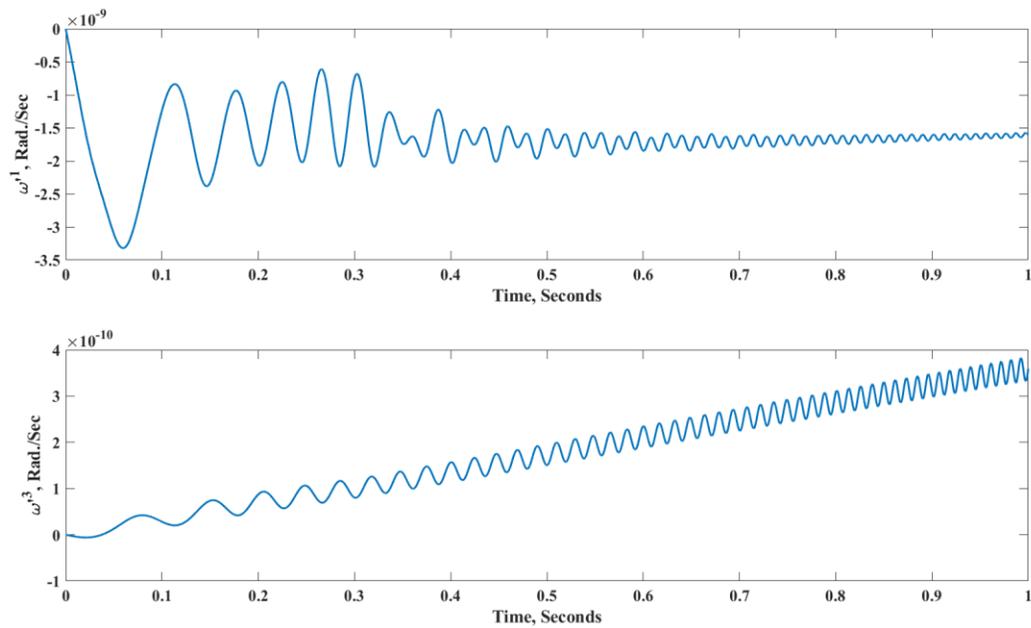

Fig. 6 Angular velocity components of flexible plate in body-fixed $X'_1$ and $X'_3$ directions

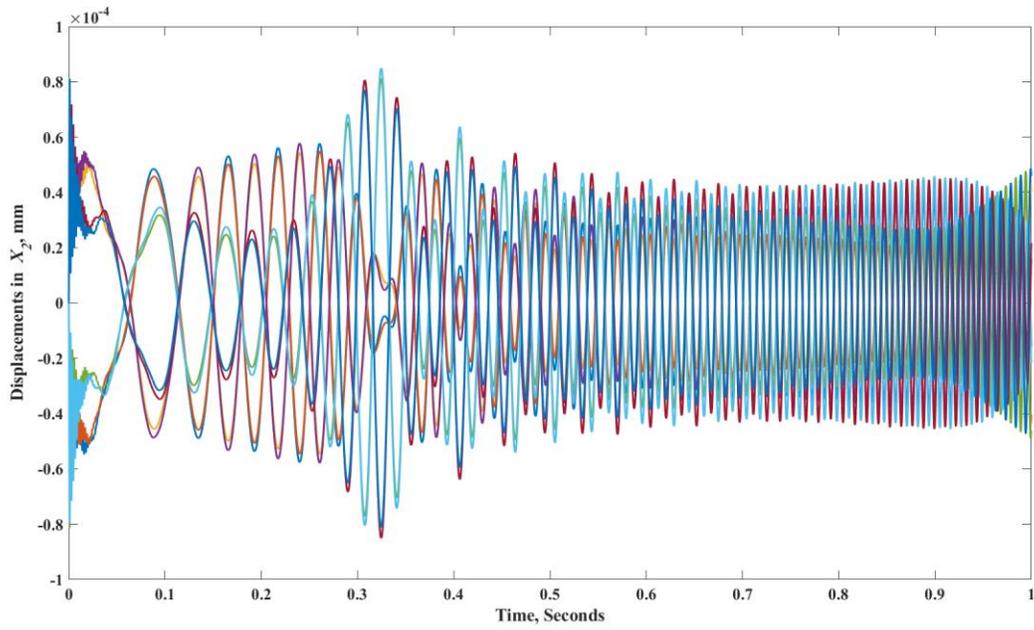

Fig. 7 Displacements of $P_1$ through $P_8$ in global $X_2$ direction

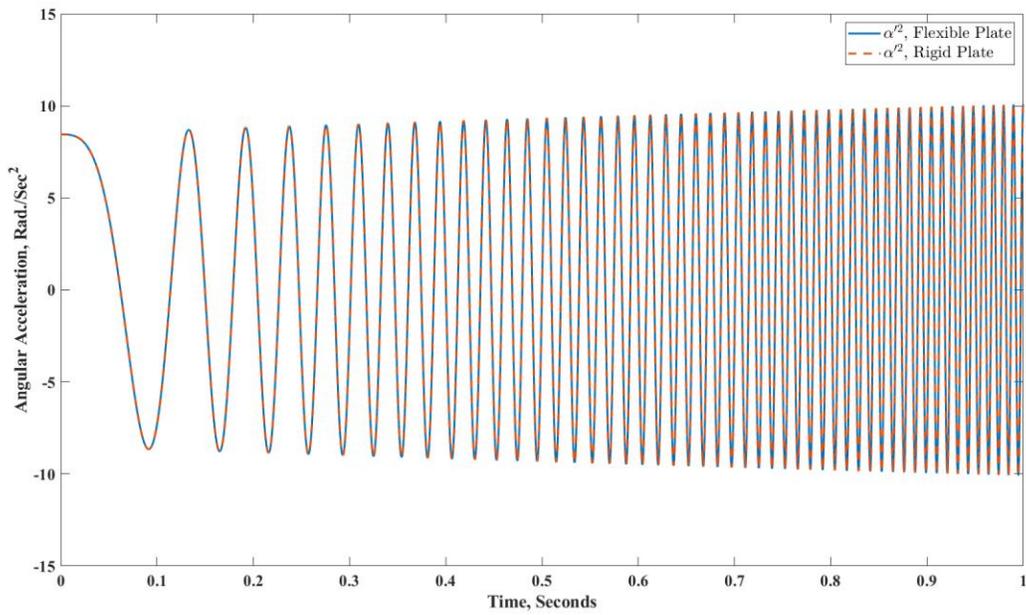

Fig. 8 Angular acceleration components in body-fixed $X_2'$ direction

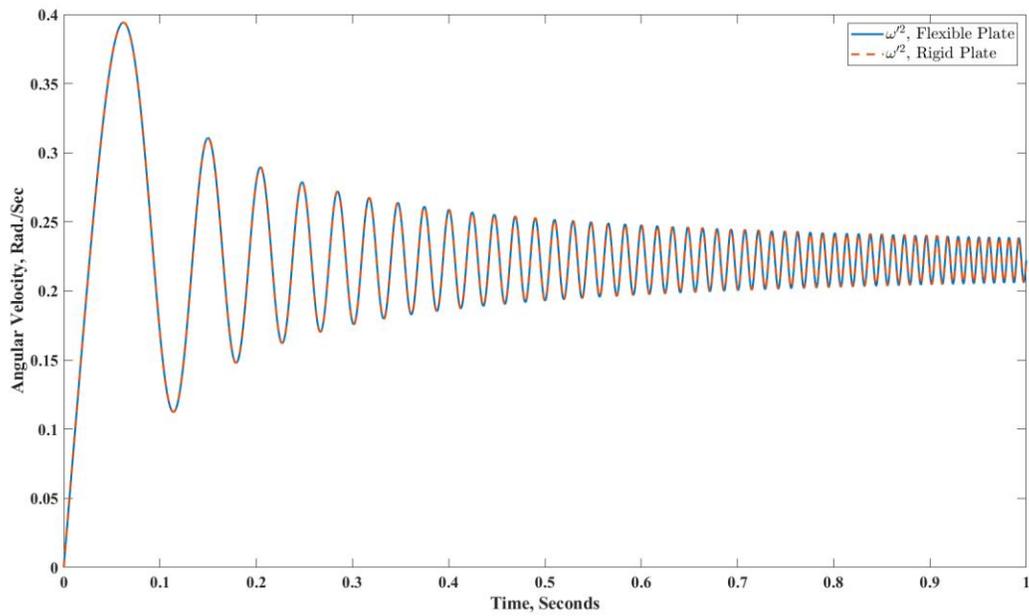

Fig. 9 Angular velocity components in body-fixed $X'_2$ direction

In Figs. 10 and 11, we compare the accelerations of point $P_1$ in global $X_1$ and $X_3$ directions obtained from the embedded RBM component of the flexible plate against the corresponding acceleration components of the rigid plate. The accelerations match precisely.

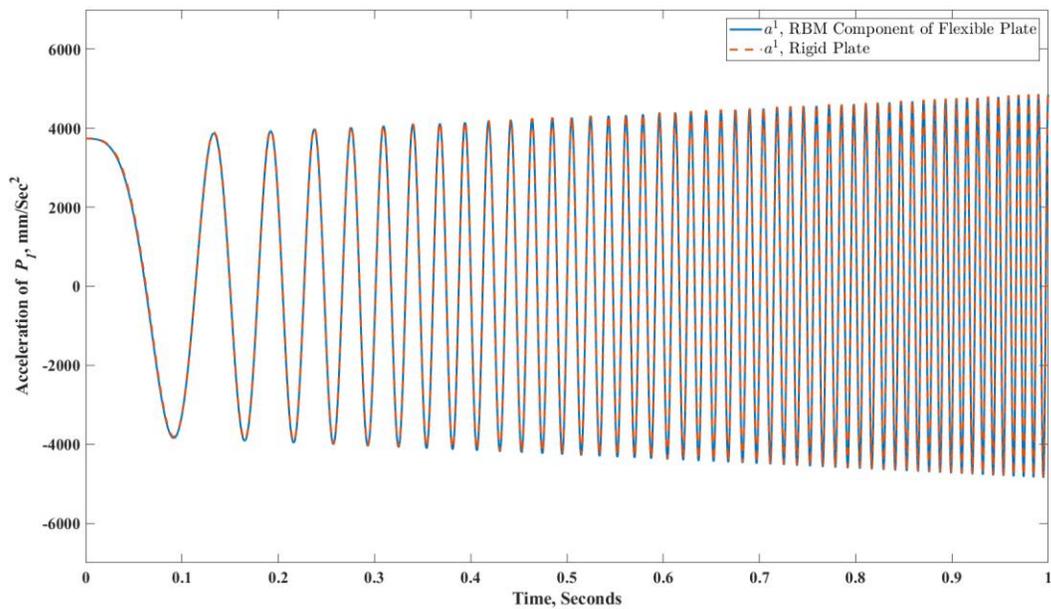

Fig. 10 Comparison of $P_1$ acceleration components in global $X_1$ direction

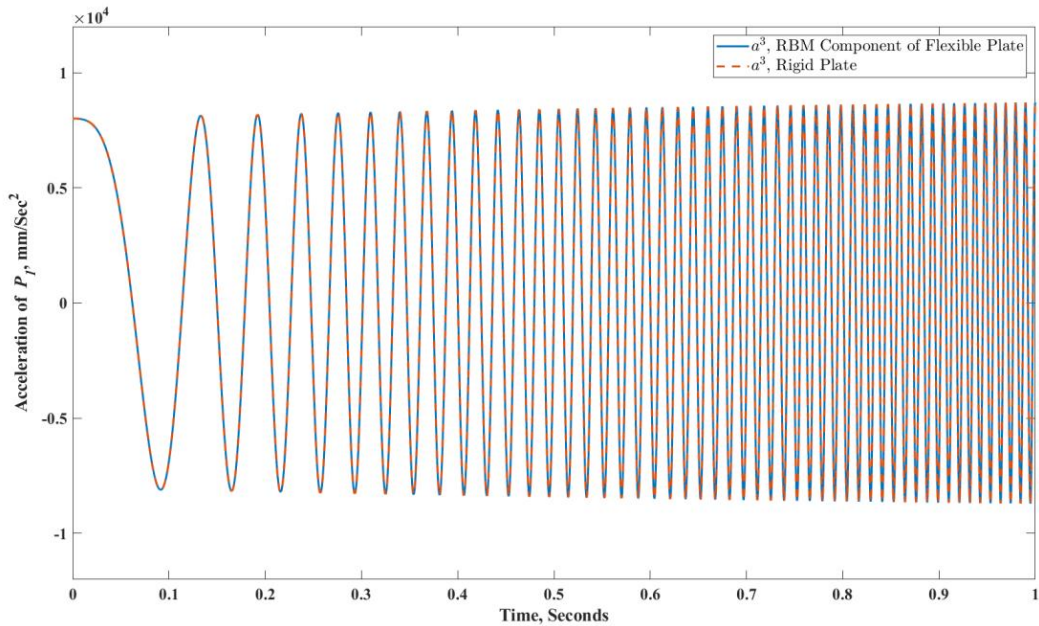

Fig. 11 Comparison of $P_1$ acceleration components in global $X_3$ direction

In Figs. 12 and 13, we compare the full acceleration components of point $P_1$ on the flexible plate in global $X_1$ and $X_3$ directions against the corresponding acceleration components from the rigid plate motion. As expected, the flexible plate accelerations demonstrate coupled rigid body motion and vibration.

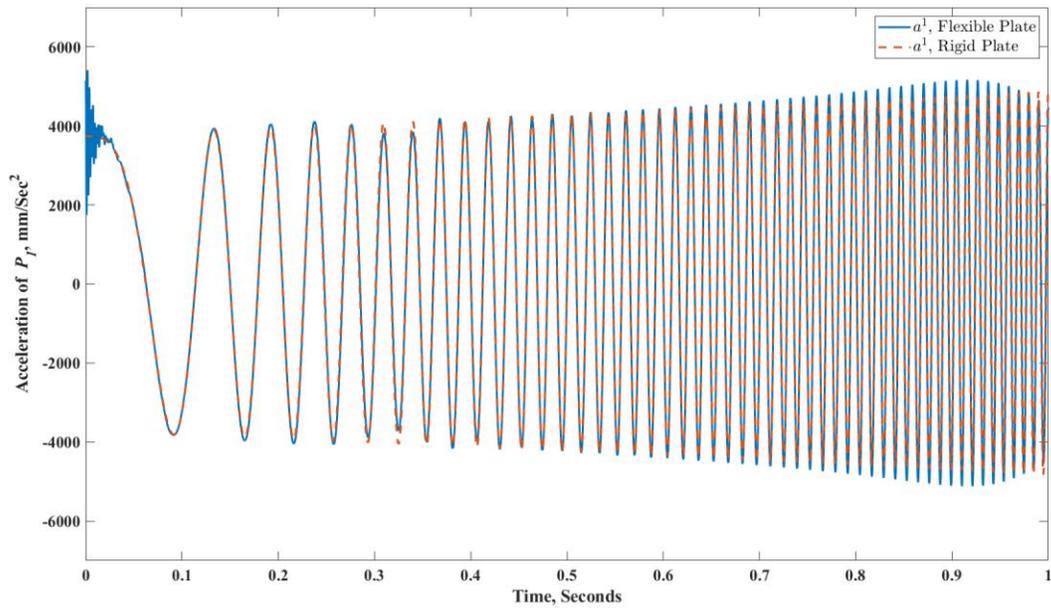

Fig. 12 Comparison of $P_1$ acceleration components in global $X_1$ direction

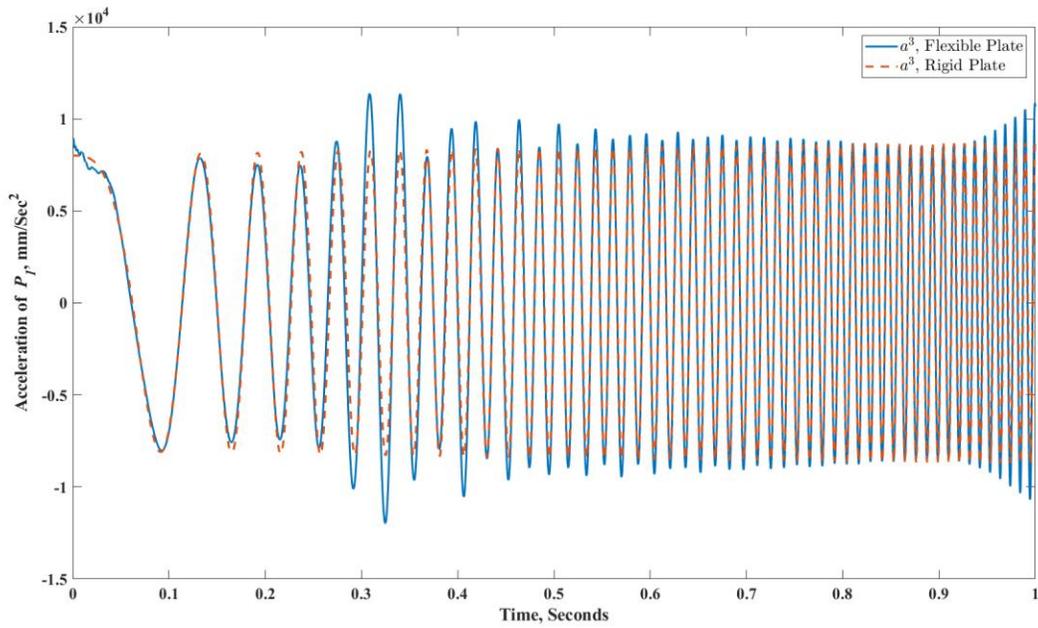

Fig. 13 Comparison of $P_1$ acceleration components in global $X_3$ direction

## 5. Conclusions

In this article we have presented a novel approach for computing flexible body dynamics. This approach extracts in an innovative way the rigid body motion component embedded within the motion of a flexible body, and it computes and synthesizes single-step responses in a sequence of reconfigured inertial frames that follow the extracted rigid body motion. The solution scheme has been rigorously validated and verified through a numerical example. Compared to established techniques such as the Floating Frame of Reference Formulation (FFRF) and the Absolute Nodal Coordinate Formulation (ANCF), this approach broadens the applicability of conventional structural dynamics models while avoiding the complexities involved in flexible body modeling and formulation of highly nonlinear equations of motion. These advantages position the new approach as a promising alternative for integration into high-performance computing tools across a wide range of application domains.

**Statements and Declarations**

This article has been authored by an employee of National Technology & Engineering Solutions of Sandia, LLC under Contract No. DE-NA0003525 with the U.S. Department of Energy (DOE). The employee owns all right, title and interest in and to the article and is solely responsible for its contents. The United States Government retains and the publisher, by accepting the article for publication, acknowledges that the United States Government retains a non-exclusive, paid-up, irrevocable, world-wide license to publish or reproduce the published form of this article or allow others to do so, for United States Government purposes. The DOE will provide public access to these results of federally sponsored research in accordance with the DOE Public Access Plan https://www.energy.gov/downloads/doe-public-access-plan.

This paper describes objective technical results and analysis. Any subjective views or opinions that might be expressed in the paper do not necessarily represent the views of the U.S. Department of Energy or the United States Government.

**Data availability** No data were used for the research described in the article.

**Conflict of interest** The author declares no competing interests.